\documentclass[12pt,thmsa]{article}
\usepackage{amsfonts}
\usepackage{amssymb}
\newtheorem{teo}{Theorem}[section]

\newtheorem{lema}[teo]{Lemma}

\newtheorem{proposition}[teo]{Proposition}

\newtheorem{obs2}[teo]{Remark}

\newtheorem{tea}{Theorem}[subsection]

\newtheorem{no2}[teo]{Note}

\newtheorem{no3}[tea]{Note}

\newcommand{\Gal}{{\rm Gal}}

\newcommand{\mod}{{\rm mod}}

\newcommand{\PSL}{{\rm PSL}}

\newcommand{\PGL}{{\rm PGL}}

\newcommand{\GL}{{\rm GL}}

\newcommand{\Image}{{\rm Image}}
\newcommand{\F}{{\mathbb{F}}}
\newcommand{\Q}{{\mathbb{Q}}}

\newcommand{\End}{{\rm End}}

\begin{document}
\title{{\bf The level $1$ weight $2$ case of Serre's conjecture
}}
\author{Luis Dieulefait\thanks{Research supported by project MTM2006-04895, MECD, Spain}
\\
Dept. d'Algebra i Geometria, Universitat de Barcelona\\
\newline
\\
{\it a Ana Chepich, ``la chica del 17"} \\ {\it (24 de
Marzo de 1917 - 8 de Julio de 2004): In memoriam}
 }
\date{\empty}

\maketitle

\vskip -20mm

\begin{abstract}  We prove Serre's conjecture for the case of Galois representations of
 Serre's weight $2$ and level $1$. We do this by combining the potential modularity
 results of Taylor and lowering the level for Hilbert modular forms with a Galois descent
  argument, properties of universal deformation rings, and the non-existence of
   $p$-adic Barsotti-Tate conductor $1$ Galois representations proved in [Di3].
\end{abstract}

\section{Introduction}

In this article we prove the non-existence of odd, two-dimensional, irreducible representations of the absolute Galois group of
$\Q$ with values in a finite field of odd characteristic $p$, in the case of Serre's weight $2$ and level (conductor) $1$. \\
 Equivalently, we prove modularity of such representations, thus solving Serre's conjecture (cf. [Se]) for them; non-existence follows from the fact that $S_2(1) = \{ 0 \}$.\\
 We will prove the result for $p>3$, the case of $p=3$ was solved by Serre (cf. [SeIII], pag. 710), based on a previous analogous non-existence result of Tate for the case $p=2$ (cf. [Tat]). In fact, in our proof we will at some step
 need the validity of the result for $p=3$: it is used to prove the modularity of certain compatible family of $\lambda$-adic conductor $1$ representations (cf. [Di3]).\\

 Before going on, let us recall to the reader the general statement of Serre's conjecture (cf. [Se]). Perhaps a good starting point is Deligne's proof of the existence of compatible families of $\lambda$-adic Galois representations attached to classical modular forms of arbitrary level and weight (cf. [De]), whose corresponding residual representations are precisely the kind of Galois representations we are interested in: odd, two-dimensional, with values on a finite field and, for almost every prime, irreducible. The idea of Serre, which appears first in print in his 1987 celebrated paper [Se], is that a converse result should be true, namely, every such two-dimensional Galois representation, with values on a finite field, should come from a classical modular form in the sense that the representation should be attached to it via Deligne's construction. This is known as the ``weak version" of Serre's conjecture. He gives also a more precise version of the conjecture, known as the ``strong version", which  specifies that a couple of Galois-theoretic invariants of a given Galois representation known as the Serre's level and Serre's weight should agree with the level and weight of at least one of the modular forms corresponding to this representation (cf. [Se]). The recipes for the predicted level and weight are clearly influenced by several ideas that were emerging at that moment, for example, by Frey's strategy regarding modularity of elliptic curves and Fermat's Last Theorem (as for the level) and by Fontaine-Laffaille's description of residual inertial weights of crystalline representations of small weight and Fontaine-Messing's comparison theorem between \'{e}tale and crystalline cohomologies (as for the weight). Serre's essential contribution is to synthesize all this information and make it part of his previously vaguer ideas regarding modularity, thus providing a more precise and much more useful conjecture. It should be emphasized that in some sense by making such an effort in giving us such a precise conjecture Serre was making an important contribution to its solution. \\
 The first spectacular result in connection to Serre's conjecture was the proof that the strong version is equivalent to the weak version. The proof of this result in full generality follows from the work of many people, the most important contributions being the control of the weight obtained by Edixhoven (see [Ed]) and the results of level-lowering of Ribet (cf. [Ri1]). It should be stressed that due to this result of Ribet the proof of Fermat's Last Theorem was reduced to the proof of the semistable case of the Taniyama-Shimura-Weil conjecture (in agreement with Frey's strategy), later proved by Wiles in [Wi]. \\
 In this paper, when proving the first cases of Serre's conjecture, we will prove directly the strong version, in fact this is inherent to our method. Nevertheless, one of the tools that we will use is a generalization of Ribet's level-lowering result to the case of Hilbert modular forms.\\
 The only cases of Serre's conjecture that were proved previous to our work are the cases already mentioned of Serre's level $1$ in characteristic $2$ and $3$ proved by Tate and Serre. A similar result was also proved for $p=5$ in [Br], assuming the Generalized Riemann Hypothesis.\\

    Results similar to ours have been obtained independently and in the same period of time by Khare and Wintenberger (cf. [K-W]).\\

 What follows is a description of the strategy of our proof: given a residual conductor $1$ representation $\hat{\rho}$ of Serre's weight $2$, we use
   results of Ramakrishna to construct a  $p$-adic representation $\rho$ deforming $\hat{\rho}$ which is
  Barsotti-Tate at $p$ and semistable at its finite set of ramifying primes. Then, the potential modularity results of Taylor imply that this representation when restricted to some totally real number field $F$ agrees with a representation attached to a Hilbert
  modular form over $F$. Over $F$, we apply the ``lowering the level" results of Jarvis, Rajaei and Fujiwara to obtain a conductor $1$
  modular $p$-adic deformation of the restriction to $F$ of $\hat{\rho}$.  Using a Galois descent argument and properties of universal deformation rings, we will see that among these conductor 1 modular $p$-adic representations there is at least one that can be extended to
    a conductor $1$ Barsotti-Tate $p$-adic deformation $\mu$ of $\hat{\rho}$. As we will see in section 3, this is an automatic consequence of our Galois descent argument and
    a result of Boeckle.  The entire process can be described as ``lowering the conductor"
    of  potentially modular Galois representations.\\
     To conclude the proof, we recall that in a previous article (see [Di3]) we have shown modularity, and therefore non-existence, of such $p$-adic representations, a result that was predicted by the Fontaine-Mazur conjecture.  For the reader's convenience, we recall the main ingredients in the proof of that result (cf. [Di3]): using potential modularity, one can build a compatible family of Barsotti-Tate conductor $1$ representations containing $\mu$, and
     looking at a prime above $3$ we derive modularity of this family from the following two facts: the residual $\mod \; 3$
      representation must be reducible, as follows from the result of Serre for $p=3$ alluded above, and then a result of
       Skinner and Wiles (the main result of [S-W]) gives modularity of one, thus all, the representations in the family.\\

       We stress that the main new idea introduced in this article is the use of potential modularity to obtain from lowering the level of Hilbert modular forms the existence of ``minimal $p$-adic deformations" of certain  residual Galois representations.
       This reduces the problem of proving cases of Serre's conjecture with small ramification to proving modularity of
        $p$-adic Galois representations of small ramification. The advantage is that in the $p$-adic case one can use the
        ``existence of  families" result proved in [Di3] which, in the case of conductor $1$ or in more general cases, can be used to switch to a small
         prime and deduce modularity there.    \\

      Remark: Except for section 3, the results in this article appeared already in a preliminary version written in March 2004 (the Galois descent argument proposed in that version was incomplete, as remarked by Kevin Buzzard in May 2004).\\
      Finally, the argument presented in section 3 to complete the proof using good properties of universal deformation rings
      dates from September 2004 (first the proof was thought to be ``conditional to a proof that the ring $R$ is not too small" because the results of Boeckle were not yet known to the author). \\

      Acknowledgements: I want to thank Xavier Xarles, Jorge Jimenez, Joan-Carles Lario, Ren{\'e} Schoof, Santiago Zarzuela, Nuria Vila, and specially Kevin Buzzard for several useful conversations. Thanks are also due to Jean-Pierre Serre for some comments on an earlier version of the preprint.\\
      I would also like to thank J.-F. Boutot, J. Tilouine and J.-P. Wintenberger, organizers of the conference ``Galois Representations" held at Strasbourg on July 2005, and J. Porras, I. Fesenko and the other organizers of the ``International Conference on Arithmetic Geometry and Number Theory" held at the Euler Institute in Saint Petersburg on June 2005, for inviting me to present these results in these conferences; and also  R. Schoof, K. Ribet and M. Dimitrov for inviting me to give a talk in the Number Theory Seminars at Rome, Berkeley and Caltech (respectively). \\
      I also want to thank the referees for their useful comments.\\

\section{Potential modularity, lowering the level and Galois descent}

Let $p>3$ be a prime, and let $\hat{\rho}$ be an odd, irreducible, two-dimensional Galois representation with Serre's weight $2$ and level $1$, with values on a finite field $\F_q$ of characteristic $p$, i.e., such that the values of the representation are $2$-by-$2$ matrices with coefficients in $\F_q$. The ``level", or ``conductor", is defined as in [Se] to be the prime-to-$p$ part of the Artin conductor, see [Se] for the definition of the weight. \\
It is well-known that such a Galois representation has a model over its ``field of coefficients", the field generated by the traces of all matrices in the image, so we will assume that the field of coefficients is $\F_q$. Observe that since the conductor is $1$ and the weight is $2$, we have $\det ( \hat{\rho}) = \chi$, the $\mod \; p $ cyclotomic character.\\
Applying the main result of [Ra] (see also [Ta4] and [K-R] for variations of this result) we know that there exists a Galois representation
$$ \rho: G_\Q \rightarrow \GL_2 (\mathcal{O}_\wp)$$
deforming $\hat{\rho}$, where $\mathcal{O}$ is the ring of integers of some number field and $\wp$ a prime of $\mathcal{O}$ above $p$, such that $\rho$ is Barsotti-Tate at $p$ (i.e., it is crystalline with Hodge-Tate weights $0$ and $1$), unramified outside
$p$ and a finite set of primes $S$ and semistable at every prime $q \in S$. Observe that in particular it holds $\det(\rho) = \chi$ (by abuse of notation, we denote by $\chi$ both the $\mod \; p$   and the $p$-adic cyclotomic characters). For such a representation $\rho$, the results of Taylor (cf. [Ta3], Theorem B or Theorem 6.1, and also [Ta2], Theorem B) imply that
there exists a totally real Galois  number field $F$ where $p$ is totally split such that the restriction of $\rho$ to $G_F$ is modular, i.e., there exists $h$ a Hilbert cuspidal modular form over $F$ of parallel weight $2$ such that
 $$ \rho|_{G_F} \cong \rho_{h, \wp'} $$
 for a prime $\wp'$ above $p$ in the field $\Q_h$ generated by the eigenvalues of $h$, where $\rho_{h, \wp'}$ denotes a representation in the family associated to $h$ in [Ta1]. \\
 So the restriction $\hat{\rho}|_{G_F}$ has conductor $1$ and it has a modular semistable deformation corresponding to a cuspform $h$ of parallel weight $2$. Thus, we can apply lowering the level results of Jarvis, Fujiwara and Rajaei for Hilbert modular forms  (see [Ja], section 1 for a description of the available lowering the level results and section 8 for an application to a semistable case, and the references therein) and conclude that there exists a level $1$ Hilbert cuspform $h'$ over $F$, also of parallel weight $2$, such that
 the corresponding Galois representation $\rho_{h', \wp''}$ gives a minimally ramified deformation of $\hat{\rho}|_{G_F}$. \\
 Thus, we have obtained a Barsotti-Tate $p$-adic conductor $1$ deformation of $\hat{\rho}|_{G_F}$. \\
 Observe that if $E \subseteq F$ is a field with $\Gal(F/E)$ solvable, using solvable base change (cf. [Ta3], Theorem 6.1) we also have potential modularity over $E$, and the above procedure gives a $p$-adic conductor $1$ modular deformation of $\hat{\rho}|_{G_E}$. In general
 such a field $E$ will not be a Galois extension of $\Q$.\\

 Another important fact, proved by Ribet (cf. [Ri2]) is that since $\hat{\rho}$ has weight $2$ and conductor $1$ its image must be
 ``as large as possible",       more precisely:
  $$ \Image (\hat{\rho}) \cong \{  x \in \GL_2(\F_q) \; : \; \det(x) \in \F_p^* \} $$
From this we obtain that the image of the projectivization of $\hat{\rho}$ is $\PGL_2(\F_q)$ or $\PSL_2(\F_q)$.\\
Since $\hat{\rho}$ is odd, the group $\PSL_2(\F_q)$ is simple and the field $F$ is totally real and $p$ is totally split in $F$, we easily see that $F$ is linearly disjoint from the field fixed by the kernel of $\hat{\rho}$. Therefore:

$$ \Image(\hat{\rho}|_{G_F}) \cong \Image(\hat{\rho}) \qquad \qquad (*) $$
 In this (and any similar) situation we will say that ``$F$ is linearly disjoint from
$\hat{\rho}$".\\

Let $S_p$ be a Sylow $p$-subgroup of $\Gal(F/\Q)$ (which may be trivial), and $E \subseteq F$ the corresponding fixed field. \\
Since $S_p$ is solvable, as we mentioned before we can also construct over $E$ a conductor $1$ Barsotti-Tate modular $p$-adic deformation of
$\hat{\rho}|_{G_E}$. \\

Therefore, without loss of generality, we will assume that we have potential modularity, and existence of minimal $p$-adic deformations,
 over a field such that $[F : \Q]$ is prime to $p$ (but $F$ may not be Galois).\\

 Since $\hat{\rho}|_{G_F}$ is absolutely irreducible, for any minimal deformation $\mu'$ of it that can be extended to a 2-dimensional representation of $G_\Q$, two different extensions will differ by a twist by a character unramified outside the primes that ramify in $F$, and in particular one and only one will correspond to a deformation of $\hat{\rho}$ (we are using the fact that
 $[F : \Q]$ is prime to $p$ and that the kernel of reduction is a pro-$p$ group, which implies that the reduction mod $p$ of the twisting character can not be trivial).\\
  Now, for such a minimal deformation $\mu'$ that extends to $G_\Q$, if we call
  $\mu$ the corresponding deformation of $\hat{\rho}$, we can see that $\mu$ is a minimal deformation: again, this follows from the fact that $[F : \Q]$ is prime to $p$ and the kernel of reduction is a pro-$p$ group, because $\mu|_{G_F} = \mu'$ has conductor $1$ and $\mu$ is a deformation of $\hat{\rho}$, also of conductor $1$. A similar argument also shows that the equality (*) of residual images implies the equality of the images of $\mu$ and $\mu'$. \\
  Observe finally that if $\mu'$ is Barsotti-Tate the same also holds for $\mu$ because $p$ is totally split in $F / \Q$.\\
  Let us record what we have proved in the following lemma:

  \begin{lema}
  \label{teo:CC}
 There is a one-to-one correspondence between minimally ramified Barsotti-Tate deformations of $\hat{\rho}$ and  those minimally ramified Barsotti-Tate deformations of $\hat{\rho}|_{G_F}$ that can be extended to the full $G_\Q$.
 Furthermore, for a deformation of $\hat{\rho}|_{G_F}$ that can be extended, the image is not enlarged in the descent process.
 \end{lema}

  In the next section, we will conclude the proof that at least one of the minimal modular Barsotti-Tate deformations of
  $\hat{\rho}|_{G_F}$ that we have obtained from lowering-the-level can be extended to $G_\Q$. This and the above lemma give:\\

Conclusion:  $\hat{\rho}$ admits a conductor $1$ Barsotti-Tate $p$-adic deformation (assuming $p>3$). As explained in section 1, we proved in a previous article (cf. [Di3]) that such representations do not exist, by proving their modularity. Thus, we conclude that
$\hat{\rho}$ can not exist. This result was proved for $p=3$ by Serre in [SeIII], pag. 710. Let us state our result in the following:

\begin{teo}
\label{teo:Serre}
For any odd prime $p$, an odd, irreducible, two-dimensional Galois representation of $G_\Q$ with values in a finite extension of $\F_p$ having conductor $1$ and Serre's weight equal to $2$ is modular, therefore it can not exist.\\
In other words: the case $k=2$ and $N=1$ of Serre's conjecture is true.
\end{teo}

\section{Minimal universal deformations and their properties}

 Let us call $R$ the universal deformation ring of minimally ramified (Barsotti-Tate at $p$) deformations of $\hat{\rho}$, and let $R'$ denote a similar minimal universal deformation ring, but of $\hat{\rho}|_{G_F}$.\\
 Let $W$ be the Witt ring of $\F_q$. We know from the results of Taylor that $R'$ is a complete intersection ring and it is
  finite flat over $W$.\\
  We will also need the following result of Boeckle (see [Bo], corollary 1, and [Ra]):\\

 \begin{proposition}
 \label{teo:gebhard} $R$ is an $W$-algebra of the type:

 $$  W [[ X_1, .... X_r]] / (f_1, ...., f_s) $$

 with $r \geq s$.
 \end{proposition}

 Now, let us explain why this finer information on minimal deformation rings is enough to conclude that in the complete intersection ring $R'$ there is at least one $p$-adic deformation that can be extended to $G_\Q$.\\
 This is straightforward because, since $R'$ is finite flat complete intersection, if we assume that none of the minimal $p$-adic deformations of $\hat{\rho}|_{G_F}$  descends to $G_\Q$, using the correspondence in lemma \ref{teo:CC} we would conclude that $R$ is too small (has Krull dimension $0$) to match with the lower bound given by proposition \ref{teo:gebhard}. \\
 This concludes the proof that some of the minimal modular $p$-adic deformations above does descend, so by lemma \ref{teo:CC} existence of minimal $p$-adic deformations, and thus Serre's conjecture in the level 1 weight 2 case, follow.\\

 Remark: We know by lemma \ref{teo:CC} that the coefficient ring of a deformation of
 $\hat{\rho}|_{G_F}$  that descends to $G_\Q$ is not increased in the descent (recall that this uses Ribet's results, only valid in the semistable weight $2$ case). \\
 From this, it easily follows from the universal property defining $R'$ that, in fact, $R$ is a quotient of $R'$. Thus, $R$ is finite, and this gives another proof of existence of minimal deformations, as proved in [Bo], lemma 2.

\section{Two Applications}
In the study of the images of certain 4-dimensional symplectic families of Galois representations, theorem \ref{teo:Serre} has important consequences. For the case of level $1$ genus $2$ Siegel cuspforms, the stronger version of the result of determination of the
 images in [Di1] was conditional to the validity of Serre's conjecture precisely in the level $1$ weight $2$ case. Thus, this result holds now unconditionally, and one can compute the images of the Galois representations attached to level $1$ Siegel cuspforms
 in any given example.\\

 In a similar way, the determination of the images of the Galois representations attached to abelian surfaces
 with $\End(A) = \mathbb{Z}$ was done in an effective way in [Di2], under the assumption of the truth of Serre's conjecture. In the case of abelian surfaces of prime conductor (there are examples constructed by Brumer and Kramer of such surfaces), it is enough to apply the case of Serre's conjecture proved in theorem \ref{teo:Serre} to determine the images unconditionally. This is due to the fact that the only case of non-maximal image requiring the validity of
Serre's conjecture (cf. [Di2]) was the case of a residually reducible representation with two two-dimensional irreducible components, both of Serre's weight $2$: if the surface has prime conductor, one of this two-dimensional components should have conductor $1$, thus contradicting our theorem.

\section{Final Comments}

If we take an odd, irreducible, two-dimensional Galois representation, with values in a finite field of odd characteristic $p$, ramified at a finite set of primes $S$ and semistable (in the sense of [Ri2]), using the ideas described in this article (assuming $p>3$) one can construct a ``minimally ramified" $p$-adic deformation, i.e., a deformation ramifying only at $p$ and $S$, also semistable at every prime in $S$. In particular, this result allows to ``lower the conductor" of a semistable potentially modular representation $\rho$, i.e., if the conductor $c$ of the residual representation $\bar{\rho}$ is a strict  divisor of the conductor of $\rho$, there exists a  semistable representation $\rho_0$ of  conductor $c$ deforming $\bar{\rho}$. \\

One can generalize the results in this paper, and the results in [Di3] of non-existence of $p$-adic representations, to prove in a similar way Serre's conjecture in the case of Serre's weight $k = 2$ and small level (conductor) $N$, with $N$ squarefree and prime to $3$: for the non-existence result in [Di3] to hold the value of the conductor must be small enough so that the result of Serre of non-existence of irreducible odd two-dimensional Galois representations with values in a finite extension of $\F_3$ and conductor $1$ extends to conductor $N$.  This extension is known to hold for example for $N=2$ and $N=7$ (cf. [Sc]). Therefore, Serre's conjecture is true for representations of weight $2$ and level (conductor) $1,2$ or $7$.\\

\section{Bibliography}

[Bo] Boeckle, G., {\it On the isomorphism $R_\emptyset \rightarrow T_\emptyset$}, appendix to: Khare, C., {\it On isomorphisms between deformation rings and Hecke rings}, Invent. Math. {\bf 154} (2003) 218-222
\newline
[Br] Brueggeman, S., {\it The nonexistence of certain Galois extensions unramified outside $5$}, J. Number Theory {\bf 75} (1999) 47-52
\newline
[De]  Deligne, P., {\it Formes modulaires et repr\'{e}sentations $\ell$-adiques}, S\'{e}m. Bourbaki 355,
 Springer Lecture Notes {\bf 179} (1971), 139-172
\newline
[Di1] Dieulefait, L., {\it On the images of the Galois representations
 attached to genus $2$
Siegel modular forms}, J. Reine Angew. Math. {\bf 553} (2002) 183-200
\newline
[Di2]Dieulefait, L., {\it Explicit Determination of the Images of the Galois Representations Attached to Abelian Surfaces with $End(A)= \mathbb{Z}$}, Exp. Math. {\bf 11} (2002) 503-512
\newline
[Di3] Dieulefait, L., {\it Existence of compatible families and new cases of the Fontaine-Mazur conjecture},
J. Reine Angew. Math. {\bf 577} (2004) 147-151
\newline
[Ed] Edixhoven, B., {\it The weight in Serre's conjectures on modular forms}, Invent. Math. {\bf 109} (1992), 563-594
\newline
[Ja] Jarvis, F., {\it Correspondences on Shimura curves and Mazur's principle at $p$}, Pacific J. Math. {\bf 213} (2004) 267-280
\newline
[K-R] Khare, C., Ramakrishna, R., {\it Finiteness of Selmer groups and deformation rings}, Invent. Math. {\bf 154} (2003) 179-198
\newline
[K-W] Khare, C., Wintenberger, J-P., {\it On Serre's reciprocity conjecture for 2-dimensional mod p representations of the Galois group of Q}, preprint, (2004); available at: http://arxiv.org/math.NT/0412076
\newline
[Ra] Ramakrishna, R., {\it Deforming Galois representations and the conjectures of Serre and Fontaine-Mazur}, Ann. of Math.
{\bf 156} (2002) 115-154
\newline
[Ri1] Ribet, K., {\it On modular representations of $\Gal(\bar{\Q} / \Q)$ arising from modular forms}, Invent. Math. {\bf 100} (1990) 431-476
\newline
[Ri2] Ribet, K., {\it Images of semistable Galois
representations}, Pacific J.  Math. {\bf 181} (1997)
\newline
[Sc] Schoof, R., {\it Abelian varieties over $\Q$ with bad reduction in one prime only}, Compositio Math. {\bf 141} (2005) 847-868
\newline
[Se] Serre, J-P., {\it Sur les repr{\'e}sentations modulaires de degr{\'e}
$2$ de $\Gal(\bar{\mathbb{Q}} / \mathbb{Q})$}, Duke Math. J. {\bf 54}
(1987) 179-230
\newline
[SeIII] Serre, J-P., {Oeuvres}, Vol. III (1986) Springer-Verlag
\newline
[S-W] Skinner, C., Wiles, A., {\it Residually reducible representations and modular forms}, Publ. Math. IHES {\bf 89} (1999) 5-126
\newline
[Tat] Tate, J., {\it The non-existence of certain Galois extensions of $\Q$ unramified outside $2$}, Arithmetic Geometry, Contemp.
 Math. {\bf 174}, A.M.S., (1994) 153-156
\newline
[Ta1] Taylor, R., {\it On Galois Representations Associated to Hilbert Modular Forms}, Inventiones Math. {\bf 98} (1989) 265-280
\newline
[Ta2] Taylor, R., {\it Remarks on a conjecture of Fontaine and Mazur},
J. Inst. Math. Jussieu {\bf 1} (2002)
\newline
[Ta3] Taylor, R., {\it On the meromorphic continuation of degree two
 L-functions}, preprint, (2001);
 available at http://abel.math.harvard.edu/$\sim$rtaylor/
 \newline
 [Ta4] Taylor, R., {\it On icosahedral Artin representations, II}, American J. Math. {\bf 125} (2003) 549-566
\newline
[Wi] Wiles, A., {\it Modular Elliptic Curves and Fermat's Last Theorem}, Ann. of Math. {\bf 141} (1995) 443-551

\end{document}